\documentclass{kluwer}
\usepackage{amsmath}
\begin{document}
\begin{opening}
\title{ MiniMax Entropy and Maximum Likelihood} 
\subtitle{Complementarity of Tasks, Identity of Solutions}
\author{Marian \surname{Grend\'ar}}
\institute{ Institute of Measurement Science, Slovak Academy of
Sciences, D\'u\-brav\-sk\'a cesta~9, 842 19 Bratislava, Slovakia
\email{umergren@savba.sk}}
\author{Mari\'an \surname{Grend\'ar}}
\institute{ Railways of Slovak Republic, DDC, Klemensova 8, 813 61 Bratislava,
Slovakia\email{grendar.marian@zsr.sk}}
\classification{AMS codes}{Primary 62-02; Secondary 62A10, 62A99, 62F10}
\abbreviations{%
\abbrev{ML}{Maximum Likelihood};
\abbrev{ME}{Maximum Entropy};
\abbrev{MiniMax Ent}{MiniMax Entropy};
\abbrev{MMM}{Modified Method of Moments};
\abbrev{FOC}{First Order Condition}
}
\begin{abstract}
Concept of exponential family is generalized by
simple and general exponential form. Simple and general potential
are introduced. Maximum Entropy and Maximum Likelihood tasks are
defined. ML task on the simple exponential  form and ME task on
the simple potentials are proved to be complementary in set-up and
identical in solutions. ML task on the general exponential form
and ME task on the general potentials are weakly complementary,
leading to the same necessary conditions. A hypothesis about
complementarity of ML and MiniMax Entropy tasks and identity of
their solutions, brought up by a special case analytical as well
as several numerical  investigations, is suggested in this case.

MiniMax Ent can be viewed as a generalization of MaxEnt for parametric
linear inverse problems, and its complementarity with ML as yet another
argument in favor of Shannon's entropy criterion.
\end{abstract}
\end{opening}

\def\vc{\mathbf}
\def\vcs{\boldsymbol}
\def\Var{\text{Var}\,}
\def\Cov{\text{Cov}\,}

\newdisplay{defn}{Definition}
\newproof{note}{Note}
\newdisplay{example}{Example}
\newproof{notation}{Notation}
\newtheorem{theorem}{Theorem}
\newtheorem{coroll}{Corollary}

\section{Introduction}

A relationship between Maximum Likelihood (ML)\index{maximum likelihood} and Maximum
Entropy \index{maximum entropy}
(ME, MaxEnt) methods has been noted and investigated many times. Yet it seems to
be intricate and puzzling. Jaynes,\cite{Jaynes82}, is worth long quoting on the
subject
\begin{quote}
..., any MaxEnt solution also defines a particular model for which
the predictive distribution using the ML estimates of the
parameters, is identical with the MaxEnt distribution. This is
essentially the Pitman-Koopman theorem used backwards; given any
data the MaxEnt distribution having exponential form, in effect
creates a model for which those data would have been sufficient
statistics. This can give one deeper understanding of the terms
'information' and 'sufficiency' in statistics, but only after some
deep thought. As  a result, almost every conceivable opinion about
the relationship between MaxEnt and ML can be found expressed in
the current literature.
\end{quote}
Some of the opinions (with different level of generality)
can be found at \cite{K}, \cite{BN}, \cite{Dutta}, \cite{Golan},
\cite{C}, \cite{N}, \cite{ali1}, \cite{ali2}. Adding to it other views
on MaxEnt itself (like interpreting Shannon's entropy function as
minus expected log-likelihood, or restrictive interpretation of
the MaxEnt recovered distribution as Maxwell-Boltzmann special
member of exponential family, or insisting on non-solvability of
Jaynes' die problem by ML method) makes investigation of
relationship between MaxEnt and ML adventurous.

In the present article we make a clear distinction between operational
mode of MaxEnt and ML methods, by defining MaxEnt task (as a simple instance of MaxEnt
method) and also ML task. An analogy between Boltzmann's deduction of
equilibrium distribution of an ideal gas in an external potential
field and probability distribution leads us to extending exponential family into {\it general exponential
form\/}, and introducing a notion of {\it simple potential} and {\it general
potential}.
Concept of complementarity is introduced, and complementarity
of ME task on simple potential  and ML task on simple exponential
form is proved. Finally, a hypothesis about complementarity of MiniMaxEnt
task on general potential and ML task on general exponential form,
suggested by a simple case analytical as well as several numerical
calculations, is put forward. The results instantaneously extends
to Relative Entropy Maximization (REM)/$I$-divergence minimization.

\section{ DEFINITIONS AND NOTATION}

The notion of exponential family is extended into simple and general exponential forms.

\begin{defn}
Let $X$ be a random variable with pmf/pdf $f_X(x)$. If $f_X(x)$
can be written in
the form of
$$
f_X(x|\vcs{\lambda}) = k(\vcs{\lambda}) e^{-U(x, \vcs{\lambda})}
$$
where $U(x, \vcs{\lambda})$ is
$$
U(x, \vcs{\lambda}) = \vcs{\lambda}'\vc{u}(x)
$$
a linear combination of functions $\vc{u}(x)$ not depending on other
parameters, and $k(\vcs{\lambda})$ is normalizing factor, then it has {\it
simple  exponential form} \index{simple exponential form}. $u(x)$ is called \index{simple potential}{\it simple potential.}

If the pmf/pdf can be written in the form of
$$
f_X(x|\vcs{\lambda}, \vcs{\alpha}) = k(\vcs{\lambda}, \vcs{\alpha}) e^{-U(x,\vcs{\lambda}, \vcs{\alpha})}
$$
where $U(x, \vcs{\lambda}, \vcs{\alpha})$ is
$$
U(x, \vcs{\lambda}, \vcs{\alpha}) = \vcs{\lambda}'\vc{u}(x,\vcs{\alpha})
$$
a linear combination of functions $\vc{u}(x,\vcs{\alpha})$ depending on
other parameters $\vcs{\alpha}$, and $k(\vcs{\lambda},
\vcs{\alpha})$ is  normalizing factor, then it has {\it general
exponential form} \index{general exponential form}. $u(x, \vcs{\alpha})$ is called
{\it general potential} \index{general potential}.

The $U(\cdot)$ function is called {\it total potential}.
\end{defn}

\begin{note}
Any class of pmf/pdf which can be written in the exponential
form is equivalently characterized by its exponential form pmf/pdf or by its
potentials.
\end{note}

\begin{example}
$\Gamma(\alpha, \beta)$ distribution has simple exponential form,
with total potential  $U(x, \vcs{\lambda}) = \lambda_1  x + \lambda_2 \ln x$;
$\lambda_1 = \frac{1}{\beta}$ and $\lambda_2 = 1 - \alpha$; $u_1(x)
= x$ and $u_2(x) = \ln x$ are the potentials. The
normalizing factor $k(\lambda_1, \lambda_2) = \frac{1}{\Gamma(1-\lambda_2)
\lambda_1^{\lambda_2 - 1}}$.

$Logistic\, (\mu, \beta)$ distribution has general exponential form with total
potential $U(x, \vcs{\lambda}, \vcs{\alpha}) = \lambda_1 u_1(x, \vcs\alpha)
+ \lambda_2 u_2(x, \vcs\alpha)$, with  $\vcs{\lambda} = [\frac{1}{\alpha_2}, 2]$,
and the potentials $u_1(\cdot) = \frac{x - \alpha_1}{\alpha_2}$,
$u_2(\cdot) = \ln(1 + e^{-\frac{x - \alpha_1}{\alpha_2}})$, and $\vcs\alpha
= [\mu, \beta]$. $k(\alpha_2) = 1/\alpha_2$.


Discrete normal distribution $dn(\lambda, \alpha)$, defined over
a support by
$$
f_X(x_i|\lambda) = \frac{e^{-\lambda (x_i - \alpha)^2}}
{\sum_i e^{-\lambda(x_i - \alpha)^2}}
$$
has total potential $U(x, \lambda, \alpha) = \lambda (x -
\alpha)^2$. It can be equivalently expressed in
simple form with $U(x, \lambda_1, \lambda_2) = \lambda_1 x + \lambda_2 x^2$,
where
$\lambda_1 = - 2\alpha\lambda$ and $\lambda_2 = \lambda$. \qquad$\diamond$
\end{example}

Standard definitions of moment and sample mean  are  extended.

\begin{defn}
$V$-moment of random  variable $X$, $\mu(V)$,  is for any function
$V(X, \vcs{\alpha})$ defined as
$$
\mu(V) = {\text E}\, V(X, \vcs{\alpha})
$$
\end{defn}

\begin{defn}
Sample $V$-moment of random variable $X$, $m(V)$, is for
any function $V(X, \vcs{\alpha})$ defined as
$$
m(V) = \sum_{i=1}^m r_i V(X_i, \vcs{\alpha})
$$
where $r_i$ is frequency of $i$-th element of support in sample.
\end{defn}

\begin{defn}
Let $\mu(V)$, $m(V)$ are $V$-moment and sample
$V$-moment, respectively. Then requirement of their equality
$$
\mu(V) = m(V)
$$
will be called {\it V-moment consistency condition}.
\end{defn}

\begin{notation}
$\vcs{\lambda}$, $\vc{u}(\cdot)$,
$\vcs{\mu}(\cdot)$ and $\vc{m}(\cdot)$ are $[J, 1]$ vectors,
indexed by $j$. $\vc x$, $\vc p$ and $\vc r$ are $[m, 1]$ vectors,
indexed by $i$, with $m$ finite or infinite. $\vcs{\alpha}$ is
$[T, 1]$ vector indexed by $t$.
\end{notation}

Since entropy maximization can be reasonably constrained by constraints
other than the moment consistency constraints (see for instance
\cite{Golan96}, \cite{M}, \cite{GJP}  or proceedings
of MaxEnt conferences), in order to be
specific, we will speak about an {\it ME task\/}. Also, {\it ML task\/} is defined.
The complementarity results obtained for the ME task
easily extends to the more general constraints used with the
Shannon's entropy maximization criterion.

\begin{defn}
{\it ML task on $f_X(x|\vcs{\theta})$.} \index{ML task}
 Let $X_1, X_2, \dots, X_n$ be a random sample from population
$f_X(x|\vcs{\theta})$.
The maximum
likelihood task on $f_X(x|\vcs{\theta})$ is to find maximum likelihood estimator
$\hat{\vcs{\theta}}$ of
$\vcs\theta$, given the sample.
\end{defn}

\begin{defn} \index{MaxEnt task}
{\it ME task on $\vc{u}(\cdot)$.} Given a sample and a
vector of known potential functions $\vc{u}(\cdot)$, the maximum
entropy task is to find the most entropic distribution $\vc p$
consistent with the set  of $\vc{u}$-moment consistency conditions.
\end{defn}

\section{ML TASK AND ME TASK}

\subsection{Simple exponential form, simple potential case}

\begin{theorem}
{\rm Complementarity of ML and ME tasks, identity of solutions\ \
} \index{complementarity}

Let $X_1, X_2, \dots, X_n$ be a random sample.
Then,

i)\ {\rm complementarity of tasks\/}

a) ML estimator  $\hat{\vcs{\lambda}}$ of $\vcs{\lambda}$
on simple exponential form $f_X(x|\vcs{\lambda}) =
k(\vcs{\lambda}) e^{- \vcs{\lambda}'\vc{u}}$
is obtained as a solution of  system of
$J$ $u_j$-moment consistency conditions,

b) the  most entropic distribution
$\vc p$  satisfying the system of $J$ $u_j$-moment consistency conditions is
the simple exponential form pmf/pdf $f_X(x| \hat{\vcs{\lambda}})$.

ii)\ {\rm identity of solutions\/}

necessary and sufficient conditions for ML task on
simple exponential form pmf/pdf $f_X(x|\vcs{\lambda}) =
k(\vcs{\lambda}) e^{- \vcs{\lambda}'\vc{u}(x)}$ and ME task on the simple potentials
$\vc{u}(x)$ are identical, and they are
$$
\mu(u_j) = m(u_j) \qquad j = 1, 2, \dots, J
$$
\end{theorem}
\begin{pf}

Discrete r.v. case.

1. ML task.
\begin{equation*}
\max_{\vcs{\lambda}}\, \,  l(\vcs{\lambda}) = \ln(k(\vcs{\lambda}))
- \sum_{j=1}^J \sum_{i=1}^m \lambda_j r_i u_j(x_i)
\end{equation*}
leads to system of $J$ first order conditions (FOC)
\begin{equation*}
\mu(u_j) = m(u_j) \qquad j = 1, 2, \dots, J
\end{equation*}

The corresponding hessian matrix of second derivatives of loglikelihood
function with respect to (wrt) $\vcs{\lambda}$ is
$$
H_{ML} = - \pmatrix
\Var(u_1) & \Cov(u_1, u_2) & \hdots & \Cov(u_1, u_J) \\
\Cov(u_2, u_1) & \Var(u_2) & \hdots & \Cov(u_2, u_J) \\
\vdots & \vdots & \ddots & \vdots \\
\Cov(u_J, u_1) & \Cov(u_J, u_2) & \hdots & \Var(u_J)
\endpmatrix
$$
negative definite,
assuring that unique global maximum was attained.

Thus, ML task on simple exponential form of pmf is identical with solving a
system of $J$ non-linear equations, the $u_j$-moment consistency conditions.

2. ME task.
\begin{equation}
\begin{split}
\max_{\vc p}\, \, H(\vc p) = &  -\sum_{i=1}^m p_i \ln p_i \\
 &\text{subject to} \\
\mu(u_j) &= m(u_j) \qquad j = 1, 2, \dots, J
\end{split}
\end{equation}
which can be accomplished by means of Lagrange\-an
$$
L(\vc p) = -\sum_{i=1}^m
p_i \ln p_i  + \sum_{j=1}^J \lambda_j (m(u_j) - \mu(u_j))
$$
leading to system of $m$ FOC
\begin{equation*}
p_i = e^{-\vcs{\lambda}'\vc{u}(x_i)} \qquad  i = 1, 2, \dots, m
\end{equation*}
which, after a normalization gives the simple exponential form as the
solution.

The corresponding hessian matrix of second de\-rivatives of the
La\-gra\-nge\-an wrt $\vc{p}$ is
$$
H_{ME} =
\pmatrix
-1/{p_1} & 0 & \hdots & 0 \\
0 & -1/{p_2} & \hdots & 0 \\
\vdots & \vdots & \ddots & \vdots \\
0 & 0 & \hdots & -1/{p_m}
\endpmatrix
$$
for $p_i > 0$ negative definite,
satisfying the sufficient conditions for a unique global maximum.

Thus, ME task with the system of $J$ $u_j$-moment consistency conditions
leads to the simple  exponential form, where  $\vcs{\lambda}$, the
Lagrange multipliers, have to be
found out of the system of nonlinear equations (1).

Continuous r.v. case.

1. ML task -- in analogy with the discrete case proof.

2. ME task.
\begin{equation*}
\begin{split}
\max_{f_X(x)}\, \, H(f_X(x)) = &  -\int f_X(x) \ln(f_X(x))\, dx \\
 & \text{subject to} \\
\mu(u_j) &= m(u_j) \qquad  j = 1, 2, \dots, J
\end{split}
\end{equation*}
which can be accomplished by means of Lagrange\-an functional
$$
L(f_X(x)) = - f_X(x) \ln(f_X(x))  + \sum_{j=1}^J \lambda_j u_j(x)
f_X(x)
$$
leading to Euler's equation (FOC)
$$
\frac{\partial L(\cdot)}{\partial f_X(x)} = 0
$$
which, after a normalization gives the simple exponential form
$$
f_X(x|\vcs{\lambda}) = \frac{e^{-\vcs{\lambda}'\vc{u}(x)}}{\int
e^{-\vcs{\lambda}'\vc{u}(x)}\, dx}
$$
\end{pf}

\begin{note}
ML task on simple exponential form and ME task on simple
potentials are complementary in the sense, that
where one starts the other one ends, and vice versa. ML starts with known
simple exponential form of pmf/pdf and ends up with ML  estimators of the
parameters, found out of the potential moment consistency equations. ME,
working on the sample,  starts with assumed form of  potential functions,
forming potential moment consistency constraints.
The most entropic distribution resolved is just the exponential form pmf/pdf ML
has assumed. And the ME estimators of its parameters are the same as the ML estimators.
We say that ML task on simple exponential form pmf/pdf and ME task on
simple potentials are {\it complementary}.

ML and ME tasks are complementary in set-up but identical in solution. Both
the tasks end up with the same mathematical problem of solving  estimators of
$\vcs{\lambda}$ out of the system of potential moment consistency equations (1).
\index{complementarity}
\end{note}

\begin{example}
Let $X_1, X_2, \dots, X_n$ be a random sample of size $n$ from discrete normal
distribution
$dn(\lambda_1, \lambda_2)$, taken in the simple exponential form.

ML task of estimation leads to solving $\lambda_1, \lambda_2$ out of
system of equations
\begin{eqnarray*}
\frac{\sum_{i=1}^m x_i e^{-(\lambda_1 x_i + \lambda_2 x_i^2)}}{\sum_{i=1}^m e^{-
(\lambda_1 x_i +
\lambda_2
x_i^2)}} &= \sum_{i=1}^m r_i x_i \\
\frac{\sum_{i=1}^m x_i^2 e^{-(\lambda_1 x_i + \lambda_2 x_i^2)}}{\sum_{i=1}^m
e^{-(\lambda_1 x_i +
\lambda_2
x_i^2)}} &= \sum_{i=1}^m r_i x_i^2
\end{eqnarray*}

which is just the system of $x$-moment and $x^2$-moment consistency
conditions.

ME task constrained by system of $x$-moment, and $x^2$-moment consistency
conditions
\begin{equation}
\begin{split}
\sum_{i=1}^m p_i x_i &=  \sum_{i=1}^m r_i x_i \\
\sum_{i=1}^m p_i x_i^2 &= \sum_{i=1}^m r_i x_i^2
\end{split}
\end{equation}
finds the most entropic distribution consistent with the constraints to have
form (after normalization)
\begin{equation}
p_i = \frac{ e^{-(\lambda_1 x_i + \lambda_2 x_i^2)}}{\sum_{i=1}^m
e^{-(\lambda_1 x_i + \lambda_2 x_i^2)}}
\end{equation}
where, $\lambda_1, \lambda_2$ should be found out of the system (2), after
plugging (3) in. \qquad$\diamond$
\end{example}

In passing we mention an identity of ML and modified method of moments (MMM) in the case of
exponential family, discovered by \cite{Huzur49} and explored
further by \cite{David74}. The identity holds also for
the simple exponential form, making ME complementary to both ML and
MMM. Note that MMM starts with a moment consistency conditions, where
understanding of moments is enhanced as done here by  Definitions  2, 3, 4.


\subsection{ General exponential form, general potential case}

Complementarity \index{complementarity} of the general exponential form  ML task and general potential ME task
can not be assessed analytically in full extent, for
sufficient conditions for maximum of likelihood or entropy function do not allow, in general,
for it.
We show, analytically, that ML task on the general exponential form  and ME task
on the general potentials  lead to the same FOC's.
This could be called 'weak complementarity'.

\begin{theorem}
Let $X_1, X_2, \dots, X_n$ be a random sample.
Then, necessary conditions for

a) ML task on  general exponential form pmf/pdf
$f_X(x|\vcs{\lambda}, \vcs{\alpha})$ \\
$= k(\vcs{\lambda}, \vcs{\alpha}) e^{-
\vcs{\lambda}'\vc{u}(x, \vcs{\alpha})}$

b) ME task on the general potentials $\vc{u}(x, \vcs{\alpha})$

are identical, and they are
\begin{equation*}
\begin{split}
\mu(u_j) &= m(u_j) \qquad j = 1, 2, \dots, J\\
\vcs{\lambda}' \vcs{\mu}\left(\frac{\partial \vc{u}}{\partial \alpha_t}\right)
 &= \vcs{\lambda}' \vc{m}\left(\frac{\partial \vc{u}}{\partial \alpha_t}\right)
\qquad  t = 1, 2, \dots T
\end{split}
\end{equation*}
\end{theorem}

\begin{pf}

Discrete r.v. case.

1. ML task.
\begin{equation*}
\max_{\vcs{\lambda}, \vcs{\alpha}}\, \,  l(\vcs{\lambda},
\vcs{\alpha}) = \ln(k(\vcs{\lambda}, \vcs{\alpha}))
- \sum_{j=1}^J \sum_{i=1}^m \lambda_j r_i u_j(x_i, \vcs{\alpha})
\end{equation*}
leads to system of $J + T$ first order conditions
\begin{equation}
\begin{split}
\mu(u_j) &= m(u_j) \qquad j = 1, 2, \dots, J\\
\vcs{\lambda}' \vcs{\mu}\left(\frac{\partial \vc{u}}{\partial \alpha_t}\right)
 &= \vcs{\lambda}' \vc{m}\left(\frac{\partial \vc{u}}{\partial \alpha_t}\right)
\qquad  t = 1, 2, \dots T
\end{split}
\end{equation}

2. ME task.
\begin{equation*}
\max_{\vc{p}(\vcs{\alpha})}\, \, H(\vc p(\vcs{\alpha})) =  -\sum_{i=1}^m p_i \ln p_i
\end{equation*}
\centerline{subject to}
\begin{equation*}
\mu(u_j) = m(u_j) \qquad j  = 1, 2, \dots, J
\end{equation*}
which can be accomplished by means of Lagrange\-an
$$
L(\vc p(\vcs{\alpha})) = -\sum_{i=1}^m
p_i \ln p_i  + \sum_{j=1}^J \lambda_j (m(u_j) - \mu(u_j))
$$
leading to system of $m + T$ FOC's
\begin{align}
p_i = e^{-\vcs{\lambda}'\vc{u}(x_i, \vcs{\alpha})}  \qquad i = 1, 2, \dots, m \qquad &\ \nonumber\\
-\sum_{i=1}^m \left(\frac{\partial p_i}{\partial \alpha_t} \ln p_i +
\frac{\partial p_i}{\partial \alpha_t}\right) + &\ \\
+  \sum_{j=1}^J \lambda_j \left(\frac{\partial m(u_j)}{\partial \alpha_t}
- \sum_{i=1}^m \left\{p_i \frac{\partial u_j(x_i,\vcs{\alpha})}{\partial
\alpha_t} +
\frac{\partial p_i}{\partial \alpha_t} u_j(x_i, \vcs{\alpha})\right\} \right)
&= 0    \qquad \forall t \nonumber
\end{align}
The most entropic distribution after normalization takes general exponential
form
$$
p_i = \frac{e^{-\vcs{\lambda}'\vc{u}(x_i, \vcs{\alpha})}}{\sum_{i=1}^m
e^{-\vcs{\lambda}'\vc{u}(x_i, \vcs{\alpha})}} \qquad i = 1, 2, \dots, m
$$
where 'ME estimators' of $\lambda$ have to be found out of the system of (5).

The $T$ of equations of the system (5)  simplifies heavily into
$$
\vcs{\lambda}' \vcs{\mu}\left(\frac{\partial \vc{u}}{\partial \alpha_t}\right)
 = \vcs{\lambda}' \vc{m}\left(\frac{\partial \vc{u}}{\partial \alpha_t}\right)
\qquad t = 1, 2, \dots, T
$$
which are  the same as the $T$ equations of FOC's  for ML task (4).

Thus, the ME and ML tasks indeed lead to the same necessary conditions (4).

Continuous r.v. case.

In analogy to the proof of Theorem 1.
\end{pf}

\begin{coroll}
Due to the linearity of $U(x, \vcs{\lambda}, \vcs{\alpha})$ in $\vcs{\lambda}$,
the necessary conditions (4) can be rewritten in a compact form
\begin{align*}
\mu\left(\frac{\partial U}{\partial \lambda_j}\right) &= m\left(\frac{\partial U}{\partial \lambda_j}\right) \qquad j = 1, 2, \dots, J\\
\mu\left(\frac{\partial U}{\partial \alpha_t}\right)
 &= m\left(\frac{\partial U}{\partial \alpha_t}\right)
\qquad  t = 1, 2, \dots T
\end{align*}
\end{coroll}

\begin{example}
Let $X_1, X_2, \dots, X_n$ be a random sample
from discrete normal distribution $dn(\lambda, \alpha)$, taken in
the general exponential form, so $u(x, \alpha) = (x - \alpha)^2$.

ML task of estimation leads to solving $\lambda$, $\alpha$ out of the
system of equations
\begin{equation}
\begin{split}
\mu(u) &= m(u) \\
\mu\left(\frac{\partial u}{\partial \alpha_t}\right) &= m\left(\frac{\partial u}{\partial \alpha_t}\right)
\end{split}
\end{equation}

ME task constrained by moment consistency condition
$$
\sum_{i=1}^m p_i (x_i - \alpha)^2 = \sum_{i=1}^m r_i (x_i - \alpha)^2
$$
leads to the FOC's
\begin{equation*}
\begin{split}
p_i &= e^{-\lambda (x_i - \alpha)^2} \\
\mu\left(\frac{\partial u}{\partial \alpha_t}\right)  &= m\left(\frac{\partial u}{\partial \alpha_t}\right)
\end{split}
\end{equation*}
where $\lambda, \alpha$ has to be found out of (6), after normalizing $p$'s.

So, ML and ME tasks lead to the same necessary conditions.
Also, note that
the ML and ME estimators are the same as in the Example 2, where $dn(\cdot)$
was taken in the simple exponential form.
\qquad$\diamond$
\end{example}

Regarding the sufficient conditions, following Theorem states  the second derivatives
for the both tasks. Whether they are identical can not be in
general analytically assessed.

\begin{theorem}
Second derivatives for the ML task are
\begin{align*}
&\frac{\partial^2 l(\vcs\lambda, \vcs\alpha)}{\partial \lambda_j^2} = - \Var(U_{\lambda_j}') \\
&\frac{\partial^2 l(\vcs\lambda, \vcs\alpha)}{\partial \lambda_j \partial\lambda_\iota} = - \Cov(U_{\lambda_j}'U_{\lambda_\iota}') \\
&\frac{\partial^2 l(\vcs\lambda, \vcs\alpha)}{\partial \alpha^2} = -\left( \Var(U_{\alpha}') + m(U_{\alpha}'') - \mu(U_{\alpha}'') \right)\\
&\frac{\partial^2 l(\vcs\lambda, \vcs\alpha)}{\partial \alpha_t\partial\alpha_\tau} =
-\mu(U_{\alpha_t}')\mu(U_{\alpha_\tau}')  - \sum_{j=1}^J \lambda_j
\mu(U_{\lambda_j\alpha_t}'' U_{\alpha_t\alpha_\tau}'') -
m(U_{\alpha_t\alpha_\tau}'') + \mu(U_{\alpha_t\alpha_\tau}'') \\
&\frac{\partial^2 l}{\partial \lambda_j\partial\alpha_t} =
-\lambda_j \Cov(U_{\lambda_j}', U_{\lambda_j \alpha_t}'') -
\sum^J_{k \ne j} \lambda_k
\mu(U_{\lambda_j}' U_{\lambda_k \alpha_t}'')
- m(U_{\lambda_j \alpha_t}'') + \mu(U_{\lambda_j \alpha_t}'')
\end{align*}

and for the ME task they are
\begin{equation*}
\begin{split}
\frac{\partial^2 L(\vc p(\vcs{\alpha}))}{\partial p_i^2} &= - \frac{1}{p_i} \\
\frac{\partial^2 L(\vc p(\vcs{\alpha}))}{\partial \alpha_t^2} &=
\Var(U_{\alpha_t}') + m(U_{\alpha_t}'') - \mu(U_{\alpha_t}'') \\
\frac{\partial^2 L(\vc p(\vcs{\alpha}))}{\partial \alpha_t\partial
\alpha_\tau} &= m(U_{\alpha_t\alpha_\tau}'') - \mu(U_{\alpha_t\alpha_\tau}'') +
\Cov(U_{\alpha_t}', U_{\alpha_\tau}')
\end{split}
\end{equation*}
\end{theorem}
\vskip3mm
\begin{pf}
Differentiating twice the loglikelihood function, and the Lagrange function lead
to the stated results.
\end{pf}

In the following simple instance of the general potential the sufficient
conditions are analytically tractable, showing  that
 at the points chosen by the necessary conditions (4) entropy
function attains its {\it maximum\/} in $\vc{p}(\vcs\alpha)$, and {\it minimum\/} in
$\vcs\alpha$,  hence the chosen distribution has minimal entropy in the
class of the most entropic distributions, consistent with the moment consistency constraints.
Likelihood function at the points  attains its maximum.

\begin{example}
Find the sufficient conditions for the Example 3 set-up.

The general total potential is $U(x, \lambda, \alpha) =  \lambda(x - \alpha)^2$,
so the potential is $u(x, \alpha) = (x - \alpha)^2$.
The second derivatives stated in the above Theorem then simplifies into
\begin{align*}
\frac{\partial^2 l(\lambda, \alpha)}{\partial \lambda^2} &= -  \Var(u) \\
\frac{\partial^2 l(\lambda, \alpha)}{\partial \alpha^2} &= -( \lambda^2\Var(u_{\alpha}') + \lambda (m(u_{\alpha}'') - \mu(u_{\alpha}'')) )\\
\frac{\partial^2 l(\lambda, \alpha)}{\partial \lambda \partial\alpha} &=
- \left(\lambda \Cov(u, u_\alpha') + m(u_\alpha') - \mu(u_\alpha') \right)
\end{align*}
for the ML task, and into
\begin{align*}
\frac{\partial^2 L(\vc p(\alpha))}{\partial p_i^2} &= - \frac{1}{p_i} \\
\frac{\partial^2 L(\vc p(\alpha))}{\partial \alpha^2} &=  \lambda^2\Var(u_\alpha') + \lambda( m(u_\alpha'') - \mu(u_{\alpha}'') ) \\
\end{align*}
for the ME task.
Furthermore, in this case
$$
m(u_\alpha'') - \mu(u_{\alpha}'') = 0
$$
and also, due to the FOC's (4)
$$
m(u_\alpha') - \mu(u_\alpha') = 0
$$
Thus, the second derivatives for the ML task form a hessian matrix
$$
H_{ML} = -  \pmatrix
\Var(u) & \lambda \Cov(u, u') \\
\lambda \Cov(u, u') & \lambda^2 \Var(u')
\endpmatrix
$$
which is negative definite, assuring in this case, that the global maximum was attained.

ME task second derivatives are
\begin{align*}
\frac{\partial^2 L(\vc p(\alpha))}{\partial p_i^2} &= - \frac{1}{p_i} \\
\frac{\partial^2 L(\vc p(\alpha))}{\partial \alpha^2} &=  4 \lambda^2\Var(x)\\
\end{align*}
showing that entropy attains its maximum in distribution $\vc p$, and minimum in $\alpha$,
at the same point where likelihood attains its maximum.

This result was also supported by numerical investigations, elucidating the behavior.
In the $\alpha$ suggested
by FOC's entropy function attains its {\it minimum\/}, whilst the maximum
is attained for an $\tilde\alpha$ degenerating $\vc p$ into an uniform distribution.
No surprise, since the value of parameter $\alpha$ of $u(x, \alpha)$ is free to choose, and
attaining the goal of maximal entropy the value is set up such that the
uniform  distribution is reached.  \qquad$\diamond$
\end{example}

The above analytically tractable case of the sufficient conditions and several
numerical investigations of more complex general potentials lead us to propose a
{\it hypothesis\/} about complementarity of ML and {\it MiniMax Entropy\/}  tasks and
identity of their solutions, under the general exponential form, general potentials.

For the sake of completeness, the MiniMax Ent task is defined.

\begin{defn}
{\it MiniMax Entropy task.} \index{MiniMaxEnt}
Given a sample and a vector of known general potentials
$\vc{u}(x, \vcs\alpha)$,  the MiniMax Entropy task is to find in the class of
all most entropic distributions $\vc p(\vcs\alpha)$ consistent with the set  of $\vc{u}$-moment consistency
conditions, a pmf/pdf with minimal entropy.
\end{defn}

\begin{note}
If the potentials are simple, MiniMax Ent task reduces into the ME task on
simple potentials.
\end{note}

\section{CONCLUSIONS}

As a way of concluding we sum up the main points of the presented
work:

1) In light of the physical analogy mentioned in the Introduction traditional statistical
notion of exponential family (see for instance \cite{Brown}, \cite{BN})
appeared to be too restrictive. An extension to {\it general exponential form\/}, driven
by the analogy was proposed. Also, simple
and general potential were introduced in the vocabulary of statistics.

2) Maximum Entropy task, as a typical instance of MaxEnt method
and Maximum Likelihood task were defined in order to make clear the difference in operational  mode of the two methods.

3) Concept of complementarity was introduced and defined (see Note 2 at the
Section 3.1). Maximum Entropy task on simple potential and Maximum Likelihood task
on simple exponential form were proved to be complementary.

4) Exploration of the complementarity of MaxEnt on general potential and ML on
general exponential form (Sect. 3.2) led to a
generalization of MaxEnt into MiniMax Ent. It was proved that
MiniMaxEnt on general potential and ML on general exponential form
lead to the same necessary conditions. Whether the conditions are
also sufficient can not be in general analytically assessed. Simple instance
of general potential (Example 4) as well as several numerical
investigations suggests that it is the case and full extent
complementarity of MiniMaxEnt on general (parametric) potential
and ML on general exponential form can be claimed.

5) Finally, we would like to note that the complemantary \index{complementarity} relationship
of MiniMaxEnt/MaxEnt task to the ML task seems to be specific property of
Shannon's entropy criterion. In \cite{my2000} it was shown, that
so-called maximum empirical likelihood (MEL) criterion constrained
by moment consistency constraints, proposed by
\cite{M} in the context of noiseless linear inverse
problem, is not complementary with ML on the MEL recovered class
of pmf/pdf.

\section{ACKNOWLEDGEMENTS}

It is a pleasure to thank
George Judge, Ali Mohammad-Djafari, Alberto Solana and Viktor
Witkovsk\'y for valuable discussions.

\end{document}